\documentclass{amsart}

\usepackage{amsmath, amsthm, amscd}
\usepackage{amsfonts,amssymb}
\usepackage{float}
\usepackage{graphics}
\usepackage{graphicx}
\usepackage[thinlines,thiklines]{easybmat}
\usepackage{multirow}
\usepackage{algorithmic}
\usepackage{algorithm}

\newcommand{\N}{\mathbb{N}} 
\newcommand{\Z}{\mathbb{Z}} 
\newcommand{\Q}{\mathbb{Q}} 
\newcommand{\R}{\mathbb{R}} 
\newcommand{\F}{\mathbb{F}} 
\newcommand{\K}{K} 
\newcommand{\C}{\mathbb{C}} 
\newcommand{\p}{\mathfrak{p}} 

\newcommand{\Cl}{\operatorname{Cl}}
\newcommand{\OO}{\mathcal{O}}
\newcommand{\OOK}{\OO_{\K}}
\newcommand{\OK}{\OO_{\K}}
\newcommand{\Nm}{\mathcal{N}}
\newcommand{\ag}{\mathfrak{a}} 

\newtheorem{theorem}{Theorem}
\newtheorem*{theorem*}{Theorem}

\renewcommand{\algorithmicrequire}{\textbf{Input:}}
\renewcommand{\algorithmicensure}{\textbf{Output:}}

\begin{document}
\title[New techniques for number field computations]{New techniques for computing the ideal class group and a system of fundamental units in number fields}     
\author{Jean-Fran\c{c}ois Biasse}\thanks{The research was carried out while both authors were working in
Sydney with the Magma group}
\address{Department of mathematics and statistics\\ 2500 University Drive NW\\ Calgary Alberta T2N 1N4}
\email{biasse@lix.polytechnique.fr}
\author{Claus Fieker}
\address{Fachbereich Mathematik\\Universit\"{a}t Kaiserslautern\\Postfach 3049\\67653 Kaiserslautern - Germany}
\email{fieker@mathematik.uni-kl.de}

\begin{abstract}
We describe a new algorithm for computing the ideal class group, the regulator and a system 
of fundamental units in number fields under the generalized Riemann hypothesis. 
We use sieving techniques adapted from the number field sieve algorithm to derive relations 
between elements of the ideal class group, 
and $p$-adic approximations to manage the loss of precision during 
the computation of units. This new algorithm is particularily efficient for number fields of small degree 
for which a speed-up of an order of magnitude is achieved with respect to the standard methods.
\end{abstract}

\maketitle 

\section{Introduction}

Let $\K=\Q(\theta)$ be a number field of degree $n$ and discriminant $\Delta$. 
In this paper, we present fast methods for computing the structure of the ideal class group 
of the maximal order $\OK$ of $K$, along with 
the regulator and a system of fundamental units of $\OO_{\K}$. 

Class group and unit group computation are two of the four principal tasks for computational algebraic number 
theory postulated by Zassenhaus (together with the computation of the ring of
integers and the Galois group). In particular, they occur in the resolution of
Diophantine equations. For example, the Pell equation 
$$T^2 - \Delta U^2 = 1,\ \ T,U\in\Z,$$
boils down to finding the fundamental unit in a real quadratic number field of discriminant $\Delta$ (see~\cite{pell_book}). 
In addition, the Sch\"{a}ffer equation
$$y^2=1^k+2^k+\hdots + (x-1)^k,\ \ k\geq 2,$$
can be solved using solutions to the Pell equation~\cite{shaffer}. Unit computations are key ingredients in solving almost all Diophantine equations, for example when solving 
Thue equations~\cite{bilu_hanrot}. On the other hand, the computation of the ideal 
class group $\Cl(\OK)$ of a number field $K$ allows in particular to provide numerical evidence 
in favor of unproven conjectures such as the heuristics of Cohen and Lenstra~\cite{lenstra_heuristic} on the ideal class group 
of a quadratic number field, Littlewood's bounds~\cite{littlewood} on $L(1,\chi)$, or Bach's bound on the minimal bound $B$ such 
that ideals of norm lower than $B$ generate the ideal class group. The class group enters also
into the computation of the Mordell-Weil group of elliptic curves with the
descent method, or the Brauer group computations for representation theory \cite{min-rep}.

In 1968, Shanks \cite{Shanks,Shanks2} proposed an algorithm relying on the baby-step giant-step method 
to compute the structure of the class number and the regulator of a quadratic number field in 
time $O\left( |\Delta |^{1/4 + \epsilon} \right)$, or $O\left( |\Delta |^{1/5 + \epsilon} \right)$ under the 
extended Riemann hypothesis \cite{LenstraShanks}. In 1985 Pohst and Zassenhaus \cite{poza}
published an algorithm that could determine the class group of arbitray number fields.
Then, a subexponential strategy for the computation of the 
group structure of the class group of an imaginary quadratic extension was described in 1989 by Hafner and McCurley \cite{hafner}. 
The expected running time of this method is bounded by $L_{\Delta}\left(\frac{1}{2} , \sqrt{2}+o(1)\right)$ where
$$L_{\Delta}(\alpha,\beta) := e^{ \beta  (\log|\Delta|)^{\alpha}( \log\log|\Delta| )^{1-\alpha} }.$$
Buchmann~\cite{Buchmann} generalized this result to the case of an arbitrary extension, the heuristic complexity being valid for 
fixed degree $n$ and $\Delta$ tending to infinity. In a recent work~\cite{biasseL13}, Biasse described an algorithm achieving 
the heuristic complexity $L_{\Delta}\left(\frac{1}{3} , O(1)\right)$ for certain classes of number fields where both the 
discriminant and the degree tend to inifinity.

In parallel of theoretical improvements, considerable efforts have been invested to make the implementations of the subexponential 
methods efficient. In the quadratic case, Jacobson~\cite{JacobsonPhd} described an algorithm based on the quadratic sieve for 
deriving relations between elements of $\Cl(\OO_{\K})$. He successfully used it for computing the class group and the fundamental 
unit of quadratic number fields. His implementation contained some of the practical improvements described in the context of 
factorization such as self initialization and the single large prime variant. This strategy was later improved by 
Biasse~\cite{biasse_chile} who used a double large prime variant and a dedicated Gaussian elimination technique. 
Attemps have been made to generalize sieving techniques to general number fields~\cite{nies}, but the proposed algorithms 
remain impractical for the sizes of discriminant of interest. 

\subsubsection*{Our contribution} In this paper, we present an algorithm based on sieving techniques 
adapted from recent implementations of the number field sieve~\cite{RSA768} for computing 
$\Cl(\OO_{\K})$ under the generalized Riemann hypothesis (GRH). We also describe a $p$-adic method for computing the regulator and a system of fundamental units. 
We show that these methods allow a significant improvement for number fields of low degree over the current state of the art.

\section{Generalities on number fields}

Let $K$ be a number field of degree $d$. It has $r_1\leq d$ real embeddings $(\sigma_i)_{i\leq r_1}$ and $2r_2$ complex 
embeddings $(\sigma_i)_{r_1 < i \leq 2r_2}$ (coming as $r_2$ pairs of conjugates). The field $K$ is isomorphic to 
$\OK\otimes\Q$ where $\OK$ denotes the ring of integers of $K$. We can embed $K$ in 
$K_\R := K\otimes \R \simeq \R^{r_1}\times \C^{r_2}, $ 
and extend the $\sigma_i$'s to $K_\R$. Let $T_2$ be the Hermitian form on $K_\R$ defined by 
$$T_2(x,x') := \sum_i \sigma_i(x)\overline{\sigma_i}(x'),$$
  and let $\| x\| := \sqrt{T_2(x,x)}$ be the corresponding $L_2$-norm. Let $(\alpha_i)_{i\leq d}$ such that 
$\OK = \oplus_i \Z\alpha_i$, then the discriminant of $K$ is given by $\Delta = \det^2(T_2(\alpha_i,\alpha_j))$. 
The norm of an element $x\in K$ is defined by $\Nm(x) = \prod_i|\sigma_i(x)|$.

To construct the ideal class group of $\OK$, we rely on a generalization of the 
notion of ideal, namely the fractional ideals of $\OK$. They can be defined as finitely generated 
$\OK$-modules of $K$. When a fractional ideal is contained in $\OK$, we refer to it as an integral ideal, which is in 
fact an ideal of $\OK$. Otherwise, for every fractional ideal $I$ of $\OK$, there exists $r\in\Z_{>0}$ such that $rI$ is integral. Fractional ideals are assumed to be non-zero.
The sum and product of two fractional ideals of $\OK$ is given by 
\begin{align*}
IJ &= \{ i_1j_1 + \cdots + i_lj_l\mid l\in \N, i_1,\cdots i_l\in I, j_1,\cdots j_l\in J\}\\
I + J &= \{ i + j\mid i\in I , j\in J\}.
\end{align*}
The fractional ideals of $\OK$ are invertible, that is for every fractional ideal $I$, there exists 
$I^{-1}:= \{ x\in K\mid xI\subseteq \OK\}$ such that $II^{-1} = \OK$.
The set of fractional ideals is equipped with a 
norm function defined by the index $\Nm(I) := [\OK:I]$ for integral ideals, 
and the natural extension $\Nm(I/J) := \Nm(I)/\Nm(J)$.
The norm of ideals 
is multiplicative and for principal ideals, it agrees with the norm
of the generator $\Nm(x\OK) = \Nm(x)$.

The ideal class group of $\OK$ is defined by 
$\Cl(\OK) := \mathcal{I}/\mathcal{P}, $
where $\mathcal{I}$ denotes the group of fractional ideals of $K$ and $\mathcal{P}\subseteq{I}$ is the subgroup of principal 
fractional ideals. We denote by $[\ag]$ the class of a fractional $\ag$ in $\Cl(\OK)$ and by $h$ the 
cardinality of $\Cl(\OK)$. Elements of $\mathcal{I}$ 
admit a unique decomposition as a power product of prime ideals of $\OK$ (with 
possibly negative exponents). An element $x\in\OK$ is said to be a unit if $(x)\OK = \OK$, or equivalently if $\Nm(x) = 1$. 
The units of $\OK$ form a multiplicative group of the form
$$U = \mu\times \langle \gamma_1\rangle \times \cdots \times \langle \gamma_r \rangle,$$
where $\mu$ is the torsion subgroup of $U$, $r:= r_1 + r_2 - 1$ and the generators $\gamma_i$ of the non-torsion part are 
called a system of fundamental units. The regulator is an invariant of $K$ which allows us to certify the calculation of 
$\Cl(\OK)$ and $U$. It is defined as $R = \operatorname{Vol}(\Gamma)$ where $\Gamma$ is the lattice generated by vectors 
of the form 
$$(c_1\log|\gamma_i|_1 , \cdots , c_{r+1}\log|\gamma_i|_{r+1}),$$
with $|x|_i := |\sigma_i(x)|$ for $i\leq r+1$, $c_1=1$ for $i\le r_1$, $c_i=2$ otherwise (for each complex embedding $\sigma_i$, if $i\leq r+1$, then 
$\bar{\sigma_i} = \sigma_j$ for some $j > r+1$).

\section{The subexponential strategy}

The idea behind the algorithm of Buchmann~\cite{Buchmann} is to find a set of ideals 
$\mathcal{B} = \{\p_1,\cdots,\p_N\}$ whose classes generate $\Cl(\OK)$, and then consider the surjective 
morphism 
$$\begin{CD}
\Z^N @>{\varphi}>> I @>{\pi}>> \Cl(\OK)\\
(e_1, \ldots, e_N) @>>> \prod_i\p_i^{e_i} @>>> \prod_i[\p_i]^{e_i}
\end{CD}$$
From the fundamental theorem of algebra, the ideal class group satisfies $\Cl(\OK)\simeq \Z^N/\ker(\pi\circ\varphi)$. Therefore, 
the knowledge of $\ker(\pi\circ\varphi)$, which has the structure of a $\Z$-lattice, enables us to derive $\Cl(\OK)$. In the 
meantime, elements of $\ker(\varphi)$ give us units as power-products of 
relations. From these units, we hope to derive a system of fundamental units of $\OK$. The subexponential strategy 
can be broken down to three essential tasks: collecting relations, calculating the class group and calculating the 
unit group. The subexponentiality is a consequence of a careful choice of $B$.


\subsection{Relation collection}

A preliminary step to the relation collection is the choice of a generating set $\mathcal{B} = \{\p_1,\cdots,\p_N\}$ of 
$\Cl(\OK)$. We choose the set of prime ideals of norm bounded by an integer $B$. The use of the Minkowski bound 
certifies the result unconditionally, but it causes the algorithm to take a time exponential in the size of $\Delta$. 
To achieve subexponentiality, many authors chose the bound of Bach~\cite{bach}, who proved that under GRH, 
$\Cl(\OOK)$ was generated by the classes of the prime ideals $\p$ satisfying
$\Nm(\p)\leq 12\log(|\Delta|)^2.$
Although asymptotically better, in practice this bound can be larger than the one described by 
Belabas et al.~\cite{belabas_bach} who stated that under GRH, the class group was generated by the classes of the 
prime ideals of norm bounded by $B$ satisfying
\begin{align*}
\sum_{(m,\p):\Nm(\p^m)\leq B}\frac{\log\Nm(\p)}{\Nm(\p^{m/2})}\left( 1- \frac{\log\Nm(\p^m)}{\log(B)}\right)  
&> \frac{1}{2}\log{|\Delta|} - 1.9n -0.785r_1\\
&+ \frac{2.468n+1.832r_1}{\log(B)},
\end{align*}
In the rest of the paper, we assume that $\mathcal{B}$ is constructed with the bound of Belabas et al. Indeed, 
Bach's bound enlarges the dimensions of the matrices that are processed during the computation of $\Cl(\OK)$, thus 
inducing a slow-down that is not compensated by the fact that the relations are found more rapidly.

During the relation collection phase, we collect relations of the form 
$$(\phi_i) = \p_1^{e_{i,1}}\hdots\p_N^{e_{i,N}},$$
where $\phi_i\in K$. We progressively build the matrix $M:= (e_{i,j})\in\Z^{k\times N}$ where $k$ is the 
number of relations collected so far. Let $\Lambda\subseteq \ker(\pi\circ\varphi)$ be the lattice generated by the 
rows of $M$. Operations on the rows of $M$ allow us to retrieve a basis for $\Lambda$ and its determinant. 
To determine if $\Lambda$ has rank $N$, we perform operations modulo a random wordsize prime $p$. In particular, 
the $LU$ decomposition of $M$ modulo $p$ allows us to identify the prime ideals that do not contribute to the 
rank of $\Lambda$. Additional relations involving these primes increase the rank of $M$, whose rows eventually 
generate a finite index sublattice of $\ker(\pi)$. To find this index, we compute the Hermite normal form (HNF) 
of $M$, that is, we perform unimodular operations encoded by $U\in\operatorname{GL}_k(\Z)$ such that 
\[ UM = \left( 
   \begin{BMAT}(@)[0.1pt,1cm,1.5cm]{c}{c.c}
   \begin{BMAT}(e)[0.1pt,0.8cm,0.8cm]{cccc}{cccc}
h_{11}& 0      & \hdots & 0      \\
\vdots & h_{22}& \ddots & \vdots \\
\vdots & \vdots & \ddots & 0      \\
*      & *      & \hdots & h_{\tiny{N}\tiny{N}}\\
  \end{BMAT} \\
\begin{BMAT}[0.1pt,1cm,0.5cm]{c}{c} 
	(0)
\end{BMAT}
\end{BMAT}
   \right),
\]
with  $\forall j < i$ : $0\leq h_{ij} < h_{jj}$ and $\forall j > i$ : $h_{ij}=0$. Once the HNF of $M$ is computed, 
adding new rows can be done very efficiently. In the meantime, the product $\prod_i h_{i,i}$ gives us an indication 
on $[\Lambda:\ker(\pi\circ\varphi)]$, as we see in \textsection~\ref{sec:unit_comp}.



\subsection{Class group computation}

Given a matrix $A\in\Z^{N\times N}$ whose rows generate $\ker(\pi\circ\varphi)$, unimodular transformations on both rows and 
columns of $A$ yield the structure of $\Cl(\OK)$. More precisely, For every non-singular matrix 
$A\in\Z^{N\times N}$, there exist unimodular matrices $U,V\in\Z^{N\times N}$ such that
$$S:= UAV = \operatorname{diag}(d_1,\hdots,d_N)$$
where $\forall i$ such that $1\leq i < N$ : $d_{i+1} | d_i$. The matrix $S$ is called the Smith normal form (SNF) of 
$A$. 



\begin{theorem}
If the rows of $A\in\Z^{N\times N}$ are a basis for $\ker(\pi\circ\varphi)$, and if $\operatorname{diag}(d_1,\hdots,d_N)$ 
is the SNF of $A$, then
$$\Cl(\OOK)\simeq \Z / d_1\Z \times \cdots\times \Z / d_N.$$  
\end{theorem}

Once enough relations have been found, the rows of $M$ generate $\ker(\pi\circ\varphi)$, and the $N$ non-zero rows of the 
HNF of $M$ are a matrix $A\in\Z^{N\times N}$ whose rows are a basis for $\ker(\pi\circ\varphi)$, and the SNF of $A$ gives us 
$\Cl(\OK)$. However, finding 
the structure of $\Cl(\OK)$ can also be done computing the SNF of a matrix which is in practice significantly smaller 
than $A$, namely the essential part of $A$. Indeed, for each matrix $H$ in HNF, there exists an index $l$ such that 
$\forall i> l,\ h_{i,i} = 1$. The upper left $l\times l$ submatrix of $H$ is called its essential part. As 
the classes of $\p_i$ for $i>l$ are generated by those of the $\p_j,\ j \leq l$, the SNF of the essential part of $A$ 
suffices to recover $\Cl(\OK)$.

\subsection{Regulator and fundamental units computation}\label{sec:unit_comp}

Computing the regulator and a system of fundamental units of $\K$ consists of 
finding kernel vectors of $M$. Indeed, if $X = (x_1,\hdots,x_{k})$ satisfies $XM = 0$, then we have
$$\left( \prod_i \phi_i^{x_i}\right) \OOK = \OK.$$
In other words, $\gamma:= \prod_i \phi_i^{x_i}$ is a unit. Every kernel vector $X$ of $M$ yields a unit, and we want to compute the group generated by all those elements as well as the regulator
of this group, defined to be zero if the group is not of full rank.
So far, finding of relations between units is mostly done using real linear algebra (LLL).
The core problem here being the numerical instability of the matrices.
This in itself is a consequence of the well-known fact that units are very large in general,
writing the fundamental unit of a real quadratic fiels explicitly with the canonical basis
needs exponentially many digits while it is always possible to find a product
representation of size polynomial in $\log |\Delta|$.
At the end of the procedure, we verify that the assumption we made on the completeness of the lattice of relations is true. To this end, we use an approximate of the Euler product
$$hR = \frac{|\mu(\K)|\sqrt{|\Delta|}}{2^{r_1}(2\pi)^{r_2}}\lim_{s\rightarrow 1} \left( (s-1)\zeta_\K (s)\right),$$
where $\zeta_\K (s) = \sum_\ag \frac{1}{\Nm(\ag)^s}$ is the usual $\zeta$-function associated to $K$. Indeed, it allows us to derive a bound $h^*$ in polynomial time under ERH that satisfies $h^* \leq  hR < 2h^*.$ If the values 
$\det(\Gamma)$ and $\det(\Lambda)$ do not satisfy this inequality, then we need to collect more 
relations.

\section{Sieving techniques}\label{sec:sieving}

In this section, we describe sieving techniques to derive relations in $\Cl(\OK)$ for general 
number fields. This is a generalization of Jacobson's results~\cite{JacobsonPhd} for quadratic 
number fields. We provide numerical data illustrating the considerable impact of these techniques 
for class group and unit group computation in the case of low degree number fields. 
Given a generating set $\mathcal{B} = \left\lbrace \p_1,\hdots\p_N\right\rbrace$ for $\Cl(\OK)$, the 
usual method for deriving relations consists of computing random exponents $\vec{e}:=(e_1,\hdots,e_N)$, $\alpha\in\OOK$ and a reduced ideal $I_{\vec{e}}$ such that 
$$\p_1^{e_1}\cdots\p_N^{e_N} = (\alpha)I_{\vec{e}}.$$
Then, every time $I_{\vec{e}}$ is $\mathcal{B}$-smooth, we obtain a relation. As the arithmetic of ideals 
is rather expensive when $n>2$, the relation search in the computer algebra software 
Pari~\cite{pari} and versions~$2.x$ for $x < 18$ of Magma~\cite{magma} consists of enumerating short elements of $I_{\vec{e}}$ via the Fincke-Pohst method~\cite{fincke_pohst}. 

Our method consists of deriving relations from smooth values of polynomials,
thus avoiding the cost of the ideal arithmetic and of the ideal reduction. Our
method for finding smooth values is based on the recent development of the
number field sieve algorithm~\cite{RSA768}. The use of trivial methods such as
trial division for finding smooth values of our polynomials would yield the
same theoretical complexity, but would be impractical for large discriminants.
The most efficient implementation of the enumeration-based strategy for finding
relations is the one of Pari. Therefore, in the following, we assess the impact
of 
our sieving method by comparing its performance with those of Pari.

\subsection{Polynomial selection}\label{sec:pol_select}

Let $\ag$ be a $\mathcal{B}$-smooth ideal of $\OOK$. In this section, we show how to provide a polynomial $P_{\ag}\in\Z[X,Y]$ of degree $n$ such that every $(x,y)\in\Z^2$ such that $P_{\ag}(x,y)$ is $B$-smooth yields a relation. Note that in theory, $\ag$ can be any ideal, however, we obtained the best results by choosing $\ag = \OK$. Let $\alpha$ and $\beta$ be two independent elements of $\ag$. Then, we create by interpolation a $P_{\ag}\in\Z[X,Y]$  such that  
$$\forall x,y\in\Z^2,\ P_{\ag}(x,y) = \Nm(x\alpha + y\beta).$$
Every time $\phi_{x,y}:= x\alpha + y\beta$ has a smooth norm, we add the relation corresponding to the principal ideal $(\phi_{x,y})$ to the relation matrix. Before applying sieving algorithms to $P_{\ag}$ to derive relations, we need to ensure that it is likely to yield enough smooth values. Polynomial selection is an important part of the number field sieve algorithm, and so it is in our algorithm. However, the specificities of our context prevent us from directly adapting the methods of NFS for selecting the sieving polynomial. First of all, we can afford to find relations with many different choices of $\alpha$ and $\beta$, whereas the choice of a sieving polynomial in the NFS algorithm is fixed. We require that our choices of $\alpha$ and $\beta$ yield polynomials with small coefficients, and that we have a sufficient randomization at the infinite places to avoid drawing $\phi_{x,y}$ spanning the same subgroup of the unit group of $\OOK$.

To randomize the choice $\alpha,\beta$, we consider random coefficients $a_1,\hdots,a_n\in\R^n$ such that $\sum_{i\leq n} a_i = 0$. For every such $n$-tuple $\vec{a}$, we define the embedding
\[   \left.  \begin{array}{cccc}
        & \ag & \longrightarrow & \R^n \\
      \psi_{\vec{a}}:  & \alpha & \longmapsto & (a_1\log|\alpha|_1 , \hdots ,  a_n\log|\alpha|_n). \end{array} \right.\]
For every choice of $\vec{a}$, the set of elements of the form  $\psi_{\vec{a}}(\alpha)$ for $\alpha\in\ag$ is a lattice $\Lambda_{\vec{a}}$ of $\R^n$ for which we can find an LLL reduced basis for the norm 
$$T^{\vec{a}}_2: (x_1,\hdots,x_n)\mapsto e^{2a_1}x_1^2 + \hdots + e^{2a_n}x_n^2.$$
For every choice of $\vec{a}$, the first two vectors $\alpha,\beta$ of an LLL reduced basis of $\Lambda_{\vec{a}}$ are potential candidates for the creation of a polynomial yielding smooth values. Every time we draw such a couple of elements of $\ag$, we need to make sure that they do not generate the same $\Z$-module as another pair previously used. To prevent this from happening, every time we draw a couple $\alpha,\beta$ by the previous method, we express them in terms of the canonical $\Z$-basis of $\OOK$. Thus, to every couple $\alpha,\beta$ corresponds the matrix $M_{\alpha,\beta}\in\Z^{2\times n}$ of their 
coordinates. 
The HNF of $M_{\alpha,\beta}$ uniquely represents the $\Z$-module spanned by $(\alpha,\beta)$. Thus, to avoid duplicates, we store a hash of the HNF of $M_{\alpha,\beta}$ in a hash table every time we use a couple $(\alpha,\beta)$ to draw relations. We summarize the procedure of the selection of a sieving polynomial in Algorithm~\ref{alg:pol_selec}

\begin{algorithm}[ht]
\caption{Polynomial selection}
\begin{algorithmic}[1]\label{alg:pol_selec}
\REQUIRE $\ag$, $(A_1,\hdots,A_n)$ , $\operatorname{HashTable}$
\ENSURE Sieving polynomial $P_{\alpha,\beta}$ corresponding to $\alpha,\beta\in\ag$
\WHILE {a new $\alpha,\beta$ has not been found}
\STATE Draw $a_1\leq A_1,\hdots,a_n\leq A_n$ at random such that $a_1+\hdots+a_n = 0$
\STATE Let $\alpha,\beta$ be the first two elements of a LLL-reduced basis of $\Lambda_{\vec{a}}$ for $\vec{a}=(a_1,\hdots,a_n)$
\STATE Compute the hash $h_{\alpha,\beta}$ of the HNF of $M_{\alpha,\beta}$
\IF {$h_{\alpha,\beta}\notin\operatorname{HashTable}$}
\STATE Compute by interpolation $P_{\alpha,\beta}\in\Z[X,Y]$ such that $P_{\alpha,\beta} ( x, y ) = \Nm(x\alpha + y\beta)$
\ENDIF
\ENDWHILE
\RETURN $\alpha,\beta,P_{\alpha,\beta}$
\end{algorithmic}
\end{algorithm}

\subsection{Line sieving}\label{sec:line_sieving}

The quadratic sieve algorithm~\cite{pomerance_qs} used to derive smooth values of a binary quadratic form generalizes 
to the case of polynomials of arbitrary degree. Its design follows from the observation that if $P\in\Z[X,Y]$ is a 
polynomial of degree $n$, then 
\begin{equation}\label{eq:line_sieve}
\forall y_0\in\Z,\ p\mid P(r_p,y_0)\Rightarrow \forall i\in\Z, \ p\mid P(r_p + ip,y_0).
\end{equation}
Given $y_0\in\Z$, we wish to find the $x\in[-I/2,I/2]$ such that $P(x,y_0)$ is $B$-smooth, where $B$ is the bound on the 
norm of the prime ideals in the factor base. Instead of trying them all, we prefer to isolate a short list of good 
candidates that we test by trial division. If $p\mid P(x,y_0)$ for many $p\leq B$, then $P(x,y_0)$ is likely 
to be $B$-smooth. From~\eqref{eq:line_sieve}, we know that once we have one root $r_p$ of $P(X,y_0)\bmod p$, then 
we can derive all the others by translation by $(p,0)$. Line sieving consists of initiating to zero an array $S$ of length 
$I$ whose cells represent the $x\in [-I/2,I/2]$. Then, for each $p\leq B$, we compute the smallest 
roots $x_p\in[-I/2,I/2]$ of $P(X,y_0)\bmod p$ and repeat 
$$S[x_p]\leftarrow S[x_p]+\log(p), \ \ x_p \leftarrow x_p + p.$$
Then, whenever $S[x]\approx \log(P(x,y_0))$ for $x\in[-I/2,I/2]$, the value $P(x,y_0)$ is likely to be $B$-smooth. We summarize 
this procedure in Algorithm~\ref{alg:line_sieve}

\begin{algorithm}[ht]
\caption{Line sieving}
\begin{algorithmic}[1]\label{alg:line_sieve}
\REQUIRE $P\in Z[X,Y]$, $I,B,y_0\in\Z$
\ENSURE Smooth values of $P(X,y_0)$ in $[-I/2,I/2]$
\STATE $L\leftarrow \varnothing$, $\forall x\in[-I/2,I/2],\  S[x]\leftarrow 0$.
\FOR {$p\leq B$}
\STATE Let $x_p$ be the smallest root of $P(X,y_0)\bmod$ in $[-I/2,I/2]$.
\WHILE {$r_p\leq I/2$}
\STATE $S[x_p]\leftarrow S[x_p]+\log(p)$, $x_p \leftarrow x_p + p.$
\ENDWHILE
\ENDFOR
\FOR {$x\in [-I/2,I/2]$}
\IF {$S[x]\approx P(x,y_0)$}
\STATE If $P(x,y_0)$ is $B$-smooth, $L\leftarrow L\cup \{x\}$.
\ENDIF
\ENDFOR
\RETURN $L$
\end{algorithmic}
\end{algorithm}

\subsection{Lattice sieving}\label{sec:lat_sieve}

Let $P_{\alpha,\beta}(X,Y)\in\Z[X,Y]$ be the sieving polynomial described in \textsection~\ref{sec:pol_select}, $B$ the bound on 
the norm of the ideals in the factor base, and $I,J\in\Z_{>0}$. 
Every couple $(x,y)\in[-I/2,I/2[\times[1,J]$ such that $P_{\alpha,\beta}(x,y)$ is $B$-smooth yields a relations. Therefore, 
one can repeat the line sieving operation on $P_{\alpha,\beta}(X,y_0)$ for every $y_0\in [1,J]$. This method is efficient when 
sieving with primes $p < I$. but when the primes are significantly larger than $I$, the root computation at Step~3 of 
Algorithm~\ref{alg:line_sieve} is often performed for nothing since there is a good chance that none of the $x\in[-I/2,I/2[$ 
will be a root of $P_{\alpha,\beta}(X,y_0)\bmod p$. A way around that is to have an array $S$ of length $IJ$ representing 
$[-I/2,I/2[$ and to fill it by line sieving methods for the primes $p < I$ and by lattice sieving for the other primes.


The lattice sieve was first described by Pollard~\cite{lattice_sieve}. Since then, it has been extensively studied and 
improved in the past 15 years, and the most recent developments of this methods yeld the factorization of RSA768 
(see~\cite{RSA768}). This strategy relies on a one-time enumeration of roots of 
$P_{\alpha,\beta}(X,Y)\bmod p$ in $[-I/2,I/2[\times[1,J]$. 
The entry $x\leq IJ$ of the array $S$ that we use to store the logarithmic contributions corresponds to the couple 
$(i,j)\in[-I/2,I/2[\times[1,J]$ where 
\begin{align*}
i &= (x - I/2) \mod I\\
j &= (x-i-I/2)/I.
\end{align*}
As in the line sieving case, every entry of $S$ is initialized to zero, and for every $p\leq B$ and every 
$(i,j)\in[-I/2,I/2[\times[1,J]$ such that $p\mid P_{\alpha,\beta}(i,j)$, we want to perform the operation 
$S[x]\leftarrow S[x] + \log p.$ 
Line sieving repeated on every line $j\leq J$ allows us to efficiently do this for $p<I$. 
%
For the others, we followed the approach of~\cite{kleinjung}, as it is done in~\cite{RSA768} for the factorization of RSA768. 
By~\cite[Prop. 1]{kleinjung}, we know that for every $p$ such that we have a root $r_p$ of $P_{\alpha,\beta}(X,1)$ modulo $p$, 
there exists a basis $\left\lbrace (a,b),(c,d)\right\rbrace $ of the couples $(i,j)$ such that $p\mid P_{\alpha,\beta}(i,j)$ 
that satisfies
\begin{itemize}
 \item $b > 0$ and $d > 0$
 \item $-I < a \leq 0\leq c < I$
 \item $c - a \geq I$.
\end{itemize}
This basis is computed via an algorithm described in~\cite{kleinjung} that relies on the continued fraction expansion of $r_p$.
%
%
To fill the array $S$, we start from $(i,j) = (0,0)$ which is a common root modulo all primes. Then, by induction, we construct the next pair $(i' , j')$ from $(i,j)$ by choosing
\begin{itemize}
 \item $(i,j) + (a,b)$ if $i\geq -a$
\item $(i,j) + (c,d )$ if $i < I - c$
\item  $(i,j) +(a,b)+ (c,d )$ if $I-c\leq i < -a$.
\end{itemize}



\subsection{Special-$q$}\label{sec:special_Q}

The sieving space $[-I/2,I/2[\times[1,J]$ only contains a limited number of couples $(i,j)$ yielding a 
smooth value. Enlarging $I$ and $J$ might cause its size to rapidly exceed the single precision. Let $q$ be a 
prime, the special-$q$ strategy consists of sieving with a polynomial $P_q$ derived from the original sieving 
polynomial $P$ such that 
\begin{align*}
&\forall (i,j)\in [-I/2,I/2[\times[1,J],\ \  \exists(x,y)\in\Z^2, P_q(i,j) = P(x,y)\\
&\forall (i,j)\in [-I/2,I/2[\times[1,J],\ \  q\mid P_q(i,j).
\end{align*}
This strategy was used by Pollard in his original paper~\cite{lattice_sieve} to sieve on the rational side, but 
most current implementations use it in the rational side as well~\cite{RSA768}. 
To create $P_q$ for a given $q$, we need a root $r_q$ of $P$ modulo $q$. Then, we find a reduced 
basis $(a_0 , b_0),(a_1,b_1)$ of the lattice spanned by the vectors $(r_q,0),(r_q , 1)$. The polynomial 
$P_q$ is then simply given by 
$$P_q(i,j) = P(ia_0 + ja_1 , ib_0 + jb_1).$$
The reduced basis is given by successive Gaussian reductions, as explained in~\cite{kleinjung}. 
Then, to sieve with a given polynomial $P$, we repeat the procedure described in~\textsection\ref{sec:lat_sieve} for 
many different polynomials of the form $P_q$. Fortunately, once the roots of $P\mod p$ for all $p\leq B$ have been 
computed, it is possible to use these values to compute the roots of $P_q\mod p$ for $p\leq B$. Indeed, 
$$P(ia_0 + ja_1 , ib_0 + jb_1)\equiv 0\mod p$$
means that there is some root $r_p$ of $P(X,1)\mod p$ such that $r_p\equiv\frac{ia_0+ ja_1}{ib_0 + + jb_1}\mod p$. This 
implies that we have $P_q ( r_p^q , 1)\equiv 0\mod p$ for 
$$r_p^q \equiv \frac{i}{j}\equiv -\frac{a_1 - r_pb_1}{a_0 - r_pb_0}\mod p,$$
which gives us a root of $P_q (X,1)\bmod p$ from $(a_0 , b_0),(a_1,b_1)$ and a root of $P(X,1)\mod p$. We summarize 
our procedure to derive relations from an ideal $\ag\subseteq \OK$ in Algorithm~\ref{alg:sieve_ideal}.

\begin{algorithm}[ht]
\caption{Sieving procedure}
\begin{algorithmic}[1]\label{alg:sieve_ideal}
\REQUIRE $\ag\subseteq\OK$, $\mathcal{B} = \{ \p\mid \Nm(\p)\leq B\}$ , $I,J\in\Z_{>0}$. 
\STATE Select $\alpha,\beta\in\OK$ and a sieving polynomial $P_{\alpha,\beta}$ with Algorithm~\ref{alg:pol_selec}.
\STATE $\forall p\leq B$, compute the roots of $P_{\alpha,\beta}(X,1)\bmod p$. 
\FOR {$q\leq B$}
\STATE Compute $P_q$ and its roots modulo the $p\leq B$ as in~\textsection\ref{sec:special_Q}.
\STATE Let $S$ be an array of size $IJ$ initialized to 0.
\FOR {$p\leq I$}
\STATE Do $S[x]\leftarrow S[x] + \log(p)$ for each $x$ representing $(i,j)\in[-I/2,I/2[\times[1,J]$ 
such that $p\mid P_q(i,j)$ by repeating Algorithm~\ref{alg:line_sieve} for each line $j\leq J$.
\ENDFOR
\FOR {$p > I$}
\STATE Calculate a basis $\left\lbrace (a,b),(c,d)\right\rbrace $ of the lattice of points in $[-I/2,I/2[\times[1,J]$ 
that are roots of $P_q(X,Y)\bmod p$ with the method of~\textsection\ref{sec:lat_sieve}. 
\STATE Do $S[x]\leftarrow S[x] + \log(p)$ for each $x$ representing $(i,j)\in[-I/2,I/2[\times[1,J]$ 
such that $p\mid P_q(i,j)$ by using the method of~\textsection\ref{sec:lat_sieve}. 
\ENDFOR
\ENDFOR
\FOR {$x\leq IJ$}
\IF {$S[x]\approx \log(P_q(i,j))$ where $x$ represent $(i,j)\in[-I/2,I/2[\times[1,J]$}
\STATE If $\log(P_q(i,j))$ is $B$-smooth, store the corresponding relation.
\ENDIF
\ENDFOR
\end{algorithmic}
\end{algorithm}


\subsection{Overall relation collection phase}

A necessary condition to compute the class group and the unit group is to produce a full-rank relation 
matrix $M$. Our sieving methods allow us to derive relations in $\Cl(\OK)$ very rapidly, but it is 
hard to force a given prime to occur in a relation. The best performance is obtained by sieving 
with the trivial ideal $\OK$. If we want to see a given prime ideal $\p\mid (p)$ occur in a relation, 
one can use the special-$q$ with $q = p$, or sieve with the ideal $\p$. However, even after using those 
methods, some prime ideals still do not contribute to the rank of $M$. Rather than sieving in random 
power-products involving missing primes, one might prefer to switch to enumeration-based methods 
to complete the relation search. To identify the primes that need to appear in a relation, we perform 
an $LU$ decomposition of the relation matrix modulo a random wordsize prime. We try to produce enough 
relations with sieving so that the rank of $M$ is 97\% of $\#\mathcal{B}$. Then we find additional 
relations with enumeration. 

\begin{algorithm}[ht]
\caption{Full rank relation matrix computation}
\begin{algorithmic}[1]\label{alg:rel_mat_cmp}
\REQUIRE $\K$, $B$
\ENSURE A full-rank relation matrix for the primes of norm bounded by $B$
\STATE $\mathcal{B}\leftarrow\left\lbrace \p\mid\Nm(\p)\leq B\right\rbrace = \{ \p_1,\cdots,\p_N\}$.
\STATE Derive $N$ relations by repeating Algorithm~\ref{alg:sieve_ideal} with $\ag=(1)$. Let $M$ be the relation matrix.
\STATE Perform an LU decomposition of $M$ and let $\operatorname{EmptyList}$ be the list of zero columns.
\FOR {$\p \in \operatorname{EmptyList}$}
\STATE Sieve with $\p$, update $M$.
\ENDFOR
\STATE Update $\operatorname{EmptyList}$ by updating the LU decomposition of $M$.
\FOR {$\p \in \operatorname{EmptyList}$}
\STATE Find a relation involving $\p$ by enumerating short elements in random power-products.
\ENDFOR
\RETURN $M$
\end{algorithmic}
\end{algorithm}

To assess the advantage of sieving over enumeration techniques, we need to isolate its contribution to 
the performances of the class group and unit group computation. To do this, we used a modified version 
of the function $\texttt{bnfinit}$ of the computer algebra software Pari that accepts in input 
a list of precomputed relations. We interfaced via SAGE this version of Pari with a developping 
version of Magma containing a function creating relations with the sieving algorithm. The Magma function 
tries to create enough relations so that the rank of $M$ be 97\% of $\#\mathcal{B}$ and passes it 
to Pari which adds new relations with enumeration methods and calculates the class group and the unit 
group. We compared the performance of this approach to the traditional $\texttt{bnfinit}$ function of Pari. 
There are two main reasons for using a hybrid version. The first one is that Pari's implementation 
of enumeration techniques is the most efficient. As these are necessary to finish the creation of the 
relation matrix after calling the sieving algorithm, it is interesting to see how the two perform together. 
Another reason for this choice is the fact that many different algorithms contribute to the computation 
of the class group and the unit group. In particular, we use time-consuming linear algebra methods  such 
as the HNF computation. Our methodology avoids the risk of seeing the influence 
of the quality of the implementation of other algorithms occuring in the class group and unit 
group computation. 

\small
\begin{table}[t]
\begin{center}
\caption{Impact of sieving on class group and unit group computation of small degree number fields}
\label{table_cmp}
 \begin{tabular}{|r|c|r|r|}
\hline
\multicolumn{1}{|c|}{$n$} &
\multicolumn{1}{c|}{$\log_2|\Delta|$} &
\multicolumn{1}{c|}{Pari} &
\multicolumn{1}{c|}{Pari+Sieving} \\
\hline
3 & 120 & 76 & 11 \\
3 & 140 & 694 & 66 \\
3 & 160 & 6828 & 333 \\
3 & 180 & 29807 & 2453 \\
\hline
4 & 120 & 38 & 7 \\
4 & 140 & 366 & 24 \\
4 & 160 & 4266 & 175 \\
4 & 180 & 31661 & 1201 \\
\hline
5 & 120 & 33 & 18 \\
5 & 140 & 295 & 64 \\
5 & 160 & 3402 & 378 \\
5 & 180 & 16048 & 2342 \\
\hline
6 & 120 & 40 & 111 \\
6 & 140 & 294 & 161 \\
6 & 160 & 1709 & 1012 \\
6 & 180 & 14549 & 8413 \\
\hline
 \end{tabular}
\end{center}
\end{table}
\normalsize

We performed our computations on a 2.6 GHz Opteron with 4GB of memory. We used a branch of 
the developping version 
2.6.0 of Pari provided by Lo\"{i}c Grenier and the developping version of Magma, interfaced via 
Sage~4.7.2. We allocated 3GB of memory to the computation made with Pari. 
For each size $d$, we drew at random 10 number fields with discriminant satisfying 
$\log_2|\Delta| = d$. For each discriminant, we computed the class group and the unit group 
with $\texttt{bnfinit}$, which we refer to as the ``Pari'' method, and with the hybrid 
version which we refer to as the ``sieving+Pari'' method. The average timings, in CPU sec 
(rounded to the nearest integer), are presented 
in Table~\ref{table_cmp}. They illustrate the impact of sieving methods for small degree number 
fields. It is very strong for degree 3,4 and 5 number field, for which we often witness a speed-up 
of a factor at least 10 while it is rather moderate for degree 6 number fields, and negligible for 
number fields of degree 7,8. Finding smooth values of a polynomial gets more difficult when we 
increase its degree, but it is not the only reason why the impact of sieving decreases with the 
degree. Indeed, for degree 6 number field, our sieving algorithm still derives relations at a competitive 
pace, but there are many linear dependencies whereas enumeration allows a more targeted search, thus 
avoiding linear dependecies. To put these improvements into perspective, we show in 
Table~\ref{table_cmp_quad} the impact of Jacobson's self initializing quadratic sieve~\cite{JacobsonPhd} 
which is implemented in Magma 2.18. The timings for "Pari" and "Pari + Sieving" are derived under the same 
setting as for Table~\ref{table_cmp}. In addition, we added the performances of of Magma 2.18 which 
uses different methods for linear algebra. Timings for the same serie of number fields 
were reported by Jacobson 
in~\cite[Table A.3]{JacobsonPhd} on a 296Mhz SUN processor (for a fair comparison one has to take into 
account the verification time since the timings of Table~\ref{table_cmp} and Table~\ref{table_cmp_quad} 
correspond to a certification 
under GRH).

\small
\begin{table}[t]
\begin{center}
\caption{Impact of the quadratic sieve on fields generated by $X^2 + 4(10^n+1)$}
\label{table_cmp_quad}
 \begin{tabular}{|c|r|r|r|r|r|r|}
\hline
\multicolumn{1}{|c|}{$n$} &
\multicolumn{1}{c|}{20} &
\multicolumn{1}{c|}{30} &
\multicolumn{1}{c|}{40} & 
\multicolumn{1}{c|}{45} & 
\multicolumn{1}{c|}{50} & 
\multicolumn{1}{c|}{54}  \\
\hline
Magma 2.18       & 0.7 &  6   &  22  &  128&   170 &  1453     \\
Pari + Sieving   & 0.5 &  5   &  44  &  271&   593 & 1085   \\
Pari             & 0.2 &  3   &  66  &  556&  1562 & 9251  \\
\hline
 \end{tabular}
\end{center}
\end{table}
\normalsize

\section{Computing the unit group}

Assume that we have created a relation matrix $(e_{i,j})$ corresponding to the relations  
$$(\phi_i) = \p_1^{e_{i,1}}\cdots \p_N^{e_{i,N}}.$$
Every kernel vector allows us to derive a unit of $\OK$. Let $\beta_1,\cdots,\beta_k$ be a generating 
system of the unit created so far. We compute a new unit $\beta'$, and we wish to find a new minimal generating 
set for $\langle \beta_1,\cdots,\beta_k,\beta\rangle$. Usually this is done by computing
(real) logarithms of the units followed by some approximate linear algebra to find a 
(tentative) relation as well as the (tentative) new basis. This then is followed by an
exact verification of the relation to guarantee correctness. The difficulty comes from the
fact that the entries in the real matrix differ vastly in size - by several orders of magnitude -
thus making it neccessary to work with a huge precision, in fact the precision is also
subexponential in the discriminant for guaranteed results.

Here, we propose to use $p$-adic logarithms instead. The key advantage comes from the
much better control of error propagation in the linear algebra: unless division by non-units
happens, linear algebra does not increase errors. However, while the correctness is based
on the unproven Leopold conjecture about the non-vanishing of the $p$-adic regulator,
this is not a problem in practice: any relation found by the $p$-adic method can easily be
verified unconditionally, thus a failure of the algorithm would provide a counter example to
Leopold's conjecture.

We start by choosing a prime $p$ such that the $p$-adic splitting field $K_p$ has moderate degree,
here we allow at most degree $2$. Then we have $n$ embeddings $\phi_i$ of $K\to K_p$, and we
define a map $L_p:K^* \to K_p^n: x \mapsto (\log \phi_i(x))_i$ where $\phi_i$ is the usual
$p$-adic logarithm extended to $K_p$. In order to estimate the necessary $p$-adic precision, 
we also need the usual real logarithmic embedding, denoted by $L:K^*\to \mathbb R^{r+1}$.
We are looking for a (rational) solution $(x_i)$ to $\sum x_i L_p(\beta_i) = L_p(\beta')$.
Using $p$-adic linear algebra we will instead get a $p$-adic solution (or a proof that $\beta'$
is independent). Using standard rational reconstruction techniques, we derive the
rational solution from the $p$-adic one and then the integral relation between the units.
In order to estimate the $p$-adic precision, we bound numerator and denominator using Cramers rule
and universal lower bounds on the logarithms of units.
The rational solution then also satisfies $\sum x_i L(\beta_i) = L(\beta')$. Let $(\alpha_i)_i$
be a basis for $\langle \beta_1, \ldots, \beta_s, \beta'\rangle$. By Cramer's rule, we write
$$x_i = \det (L_p(\beta_1), \ldots, L_p(\beta'), \ldots, L_p(\beta_s))/
  \det (L_p(\beta_1), \ldots, L_p(\beta_s))$$
Since the (unknown) $(\alpha_i)$ form a basis, we see 
that 
$$\det (L_p(\beta_1), \ldots, L_p(\beta'), \ldots, L_p(\beta_s)) / \det (L_p(\alpha_1), \ldots, L_p(\alpha_s))$$ 
is an integer and the same is true for
$L$ instead of $L_p$, thus we can write $x_i$ as a quotient of
integers. In either case, to make sense of the determinants, we will have to 
select an apropriate number of rows to make the matrices square. To bound the integers, we make
use of the Hadramat bound for $\det (L(\beta_1), \ldots, L(\beta'), \ldots, L(\beta_s))$ and
some universal lower bound for $\det (L(\alpha_i))_i$. For the lower bound we use lower bounds of
logarithms of non-torsion units (\cite[3.5]{pohst}): $\|L(\alpha_i)\|_2\ge \frac{21}{128}\frac {\log d}{d^2}$, or, if the unit
group has full rank, $s=r =r_1+r_2-1$, we use lower regulator bounds, possibly coming
from the Euler product. Having obtained bounds from the real logarithm ($L$) with
low precision, we calculate the $p$-adic precision required to find $x_i$ using
$p$-adic linear algebra and rational reconstruction. In the course of the computation
it can happen that the $p$-adic determinants ($p$-adic regulators) have non-trivial 
valuation. In this case we have to restart the computation with a correspondingly
higher precision to account for the loss.
Since the Leopold-conjecture is non-proven as of now, we also need to verify the solution
by computing a low-precision estimate for $\|\sum x_i L(\beta_i) - L(\beta')\|$ to
compare it to the lower bound used above.

From the relation $x_i$ we can easily obtain a presentation of the new basis $\alpha_i$
in terms of the $\beta_i$, $\beta'$. For optimization, we then proceed to compute
a new basis $\tilde\alpha_i$ such that the real logarithms are (roughly) $LLL$-reduced.
We note that we do not rely on any $LLL$ estimates here, so any heuristic algorithm
that aims at reducing the aparent size will do. Since we do not have any $LLL$ algorithm
that will accept real input (as opposed to rational), it is important that this 
does not influence the correctness.

\subsection{Advantages of the $p$-adic method}
There are two core advantages of the $p$-adic logarithms over the ordinary, complex, ones:
first, the linear algebra problems we need to solve in order to find dependencies or relations
between units have a much simpler error analysis. In fact, contrary to the complex case, it is
possible through the use of ring based operations to solve linear equations without any additional
loss of precision. This is very important in the context of unit computation since
the matrices representing the image of $L(\alpha)$ are very badly conditioned for
classical numerical methods. The other advantage of the $p$-adic logarithms is more subtle: if
we assume Leopold's conjecture to hold for the field(s) we are interested in, then instead
of doing linear algebra over $\R$ with a precision of say $q$ to find dependencies, it
is sufficient to work with a real precision of $q/2$ and a $p$-adic precision of $q/2$ as well.
Thus, assuming classical multiplication, we gain a factor of about 4 through the use of lower
precision. Using fast multiplication (in high precision), the gain is smaller but still
noticeable. But the most important advantage is the much easier precision control: instead of
complicated and very delicate estimates for linear algebra problems, all we need are upper bounds
on linear combinations with integral coefficients - which are trivial to obtain.

We should also mention that one disadvantage of the $p$-adic method lies in the total lack
of control over the real size of the units, thus it needs to be paired with a crude
(and uncritical for correctness) size reduction algorithm.
Also, it is (currently) not possible to forego completely the use of complex (or real)
logarithms, as the $p$-adic method is not capable to proving a unit to be torsion - without
the knowledge of bounds on the real size.

\subsection{Lower bound from Euler Product}
Suppose that, as in the class group algorithm, we are given an approximation
of the Euler product, ie. we have a real number $E$ such that $1/\sqrt 2\le hR/E\leq \sqrt 2$.
After the relation matrix has full rank, and assuming the factorbasis is large enough
for correctness, we have an upper bound for the class number, thus a lower bound for $R$.
This lower bound will be several orders of magnitude larger than the universal bounds
available otherwise.

\subsection{Saturation}
After the inital steps of the algorithm, when the relation matrix has full rank, 
we have a tentative class number $h$ and a tentative 
regulator $R$. Experimentally, at this point, $hR$ does not approximate the Euler product very well - 
the product will be off by several orders of magnitude.
However, after we found one or two more relations, the product has the same size than the
Euler product, it frequently even looks like only a factor of 2 is missing in either $h$ or $R$.
To find the last missing relation can easily take more time than the entire previous run therefore
we suggest using saturation methods instead. At this point in the algorithm
the relations
define a subgroup $U$ of the $S$-unit group $U_S$ where $S$ is the factor basis. From
the Euler product we know that the index $(U_S : U)=: b$ is small, lets say $b<B$.
For any prime $p|b$ there is some $u\in U_S\setminus U$ such that $u^p\in U$.
Let us fix the prime $p$. For any prime ideal $Q\notin S$ such that $p|\Nm(Q)-1$ we can define the map
$\phi_Q:U\to \F_Q^*/(\F^*_Q)^p$ mapping $S$-units into the multiplicative group of the
residue class field modulo $p$-th powers. The Cebotaref theorem \cite{tscheb} guarantees that if
$u\in U$ is not a $p$-th power, there will be some $Q$ such that $\phi_Q(u)$ is non-trivial, 
ie. $u$ is not a $p$-th power modulo $Q$. We now simply intersect $\ker\phi_Q$ for several $Q$
until either the intersection is $U^p$ or it is does not change for five consecutive $Q$.
We expect that any $u\in U/\cap \ker\phi_Q$ will have a $p$-th root in $U_S$ but not in $S$.
Therefore $v^p=u$ is a new relation that will change $hR$ by $p$. Repeating this
for all $p<B$ until we cannot enlarge $U$ any more we find the missing relations.

\subsection{Representation}
During the execution of the algorithm, all ($S$-)units are naturally represented
as power products of the relations coming from the sieving (or the saturation). 
It is well known that the explicit representation of the units with respect to a fixed basis 
for the field can require exponentially large coefficients, so it is important to operate
on the power products as much as possible. However, even the exponent vectors
constructed for the basis of the unit group, or the saturation, will become huge, so we need
to ``size reduce'' the power products. In particular, this happens even if the
resulting element is not too large. Using ideas of \cite{buchmann:rep} for compact representations
and \cite{hess} for reduced divisors in function fields, we can find a representation
for those elements that depends only on the logarithmic size (and the number field) rather than
the execution path. For any prime $p$ we can write any unit
$$u = \prod r_i^{e_i} = \prod a_i^{p^i}$$
with elements such that the size of $a_i$ depends on the discriminant and $p$ only. The length of 
the product comes from $L(u)$. Furthermore, in this presentation it is easy to test for
$p$-th powers as only $a_0$ needs to be tested and this is a small element.





\subsection{Example}
To illustrate the power of the $p$-adic method, we look at a totally real
quartic field generated by a root of
$$x^4 + 17211x^3 + 5213x^2 - 176910463x - 4958.$$
The discriminant $\Delta$ of the maximal order has 38 digits. In the
course of the computation, we found 534 relations involving prime
ideals of norm up to $3000 = 0.4 \log^2 |\Delta|$ describing a trivial
class group. We then searched for $5$ further relations to obtain units
$u_i$ ($1\leq i\leq 5$). As power products of the relations, the units
are given via exponent vectors $e_i$ with $\|e_i\|_\infty$ ranging between
$10^{80}$ and $10^{160}$ and $20< \|e_i\|_1/\|e_i\|_\infty < 92$. So,
while not uniformely large, the exponents are non-sparse, involving huge
integers. Using a decimal precision of 170 digits, we establish that the
logarithms of the units are roughly $\|L(u_i)\|_\infty \approx 10^{160}$.
The first three units are indeed independent, giving a basis for a
subgroup of full rank, the fourth is then dependent. Choosing the prime
$p=10337$ we get $\Q_p$ as a splitting field. Using 
a $p$-adic precision of 245 digits (ie. working in $\Z_p \bmod p^{245}$),
we compute the dependency for the fourth unit, involving exponents
of around $10^{360}$. The new unit group is then tentatively LLL reduced,
producing a new basis where the $\|L(\tilde u_i)\|_\infty$ are bounded by $10^7$
only.
The last unit then involves a much smaller dependency, here the exponents
are only around $10^{60}$.

Unfortunately, looking at the Euler product, the unit group is not
complete. However, the saturation technique outlined above takes
1sec to determine that the product of the three basis elements is
(proabably) a square. Finding a better representation where the 
exponents are all powers of 2 takes less than $1$ sec and then we can enlarge the
unit group easily.

Due to the implementation, the $p$-adic precision used was actually higher: 
changing (increasing) precision is very compuationally expensive, so
we try to avoid this and simply double the precision. We used a precision
of $320$ for the $p$-adics and a maximal precision of $1000$ for the real
precision. The computation of the $\log$ is the dominating part, we spent
50sec or 90\% of the total processing time here.





\section{Conclusion}

We introduced new techniques to enhance the performances of the subexponential methods 
for computing the class group and the unit group of a number field. In particular, sieving 
allows a speed-up of an order of magnitude for number fields of small degree. These techniques 
could be developped even further. Indeed, we have not taken into account all the improvements 
to sieving techniques describes in the context of the number fields sieve algorithm such 
as large prime variations or cache-friendly methods. It is also notable that fast techniques 
for deriving relations in the class group of a small degree number field have applications 
in evaluating isogenies between small genus curves via complex multiplication methods. Indeed, 
in that case, evaluating isogenies between genus $g$ curves involves relations in the class group 
of a degree $2g$ number field. 

\section*{Acknowledgments}

The first author is particularily grateful to Lo\"{i}c Grenier for providing a special branch 
of Pari allowing to start the class group computation from an existing relation matrix. He also 
thanks David Roe for helping with the SAGE interface between Magma and Pari.

\bibliographystyle{plain}
\bibliography{article}

\end{document}